# CONVERGENCE OF FUNCTIONALS OF SUMS OF R.V.S TO LOCAL TIMES OF FRACTIONAL STABLE MOTIONS

By P. Jeganathan

*Indian Statistical Institute*

Consider a sequence $X_k = \sum_{j=0}^{\infty} c_j \xi_{k-j}$, $k \geq 1$, where $c_j$, $j \geq 0$, is a sequence of constants and $\xi_j$, $-\infty < j < \infty$, is a sequence of independent identically distributed (i.i.d.) random variables (r.v.s) belonging to the domain of attraction of a strictly stable law with index $0 < \alpha \leq 2$. Let $S_k = \sum_{j=1}^{k} X_j$. Under suitable conditions on the constants $c_j$ it is known that for a suitable normalizing constant $\gamma_n$, the partial sum process $\gamma_n^{-1} S_{[nt]}$ converges in distribution to a linear fractional stable motion (indexed by $\alpha$ and $H$, $0 < H < 1$). A fractional ARIMA process with possibly heavy tailed innovations is a special case of the process $X_k$. In this paper it is established that the process $n^{-1}\beta_n \sum_{k=1}^{[nt]} f(\beta_n(\gamma_n^{-1} S_k + x))$ converges in distribution to $(\int_{-\infty}^{\infty} f(y)\,dy) L(t, -x)$, where $L(t, x)$ is the local time of the linear fractional stable motion, for a wide class of functions $f(y)$ that includes the indicator functions of bounded intervals of the real line. Here $\beta_n \to \infty$ such that $n^{-1}\beta_n \to 0$. The only further condition that is assumed on the distribution of $\xi_1$ is that either it satisfies the Cramér's condition or has a nonzero absolutely continuous component. The results have motivation in large sample inference for certain nonlinear time series models.

**1. Introduction.** Consider a sequence $\xi_j$, $-\infty < j < \infty$, of independent identically distributed (i.i.d.) random variables (r.v.s) belonging to the domain of attraction of a strictly stable law with index $0 < \alpha \leq 2$. Define $X_k = \sum_{j=0}^{\infty} c_j \xi_{k-j}$, where $c_j, j \geq 0$, is a sequence of real numbers. Let $S_k = \sum_{j=1}^{k} X_j$, $k \geq 1$. Then under suitable conditions on the constants $c_j$ it is known that for a suitable $0 < H < 1$ and for a suitable slowly varying function $u(n)$, the finite-dimensional distributions of the process $(n^H u(n))^{-1} S_{[nt]}$ converge in distribution to those of a linear fractional stable motion (LFSM). See, for example, Kasahara and Maejima ([1988](#)). When $\alpha = 2$, the LFSM reduces to the fractional Brownian motion (FBM), and when $H = 1/\alpha$, it is taken to be the $\alpha$-stable Lévy motion. (Definitions of the preceding processes are recalled in Section 2.)

Now, for simplicity, let $\gamma_n = n^H u(n)$. Let the constants $\beta_n$ be such that $\beta_n \to \infty$ with $n^{-1}\beta_n \to 0$. In this paper it is established that the finite-dimensional distributions of the process $n^{-1}\beta_n \sum_{k=1}^{[nt]} f(\beta_n \gamma_n^{-1} S_k + x))$ (indexed by $t$ and $x$) converge in distribution to those of $(\int_{-\infty}^{\infty} f(y)\,dy) L(t, -x)$, where $L(t, x)$ is the local time of the LFSM, for a wide class of functions $f(y)$ that, in particular, includes the indicator functions of bounded intervals. (The only further condition that will be assumed on the distribution of $\xi_1$ is that either it satisfies the Cramér's condition or has a nonzero absolutely continuous component.)

For the particular situation where the limit of $\gamma_n^{-1} S_{[nt]}$ is a Brownian motion or a FBM, some partial results in some form are available in Akonom ([1993](#)), Park and Phillips ([1999](#), [2001](#)) and Tyurin and Phillips ([1999](#)), where the motivation is an interesting development of a large sample theory in some time series models that have functions of the form $f(S_k)$ occurring as regressions. The present paper has the same motivation.

Major works for the i.i.d. situation $S_k = \sum_{j=1}^{k} \xi_j$ that are related to the approach of the present paper include Skorokhod and Slobodenjuk ([1970](#)) and Ibragimov ([1984](#), [1985](#)). The approach of Ibragimov ([1984](#), [1985](#)) [a thorough account of which is presented in the book by Borodin and Ibragimov ([1995](#))] does not rely on the moment conditions and hence, implicitly on the restriction of the Brownian motion limit, of the earlier approach of Skorokhod and Slobodenjuk ([1970](#)) but mainly on the requirement of the attraction of $\sum_{j=1}^{k} \xi_j$ to stable laws. This approach is based on the representation

$$\frac{\beta_n}{n} \sum_{k=1}^{n} f(\beta_n \gamma_n^{-1} S_k) = \frac{1}{n} \sum_{k=1}^{n} f_n(\gamma_n^{-1} S_k) = \int_{-\infty}^{\infty} \left( \frac{1}{n} \sum_{k=1}^{n} e^{-iu\gamma_n^{-1}S_k} \right) \widehat{f_n}(u)\,du,$$

where $f_n(y) = \beta_n f(\beta_n y)$ and $\widehat{f_n}(u)$ is the Fourier transform of $f_n(y)$. The conditions sought (in the i.i.d. case) are naturally through the Fourier transform $\widehat{f_n}(u)$ which in certain situations are then transformed in terms of $f_n(y)$ and/or on the distribution of $\xi_1$.

The approach of the present paper involves the approximation that the difference

$$(1) \quad \frac{1}{n} \sum_{k=1}^{n} f_n(\gamma_n^{-1} S_k) - \frac{1}{n} \sum_{k=1}^{n} \int_{-\infty}^{\infty} f_n(\gamma_n^{-1} S_k + z\varepsilon)\phi(z)\,dz \to 0 \qquad \text{in } \mathbf{L}^2$$

as $n \to \infty$ first and then $\varepsilon \to 0$, where $\phi(z) = \frac{1}{\sqrt{2\pi}} \exp\{-\frac{z^2}{2}\}$. The approximating quantity $\frac{1}{n} \sum_{k=1}^{n} \int_{-\infty}^{\infty} f_n(\gamma_n^{-1} S_k + z\varepsilon)\phi(z)\,dz$ can be handled relatively easily. This approach has some advantages, especially for the situation of the present paper, for instance, the required conditions for establishing the approximation (1) can be viewed directly through $f_n(y)$ and the distribution of $\xi_1$ themselves.

In Section 2 we recall the definition of the LFSM, state a result on the existence of its local time and also recall a result on the weak convergence of the process $\gamma_n^{-1} S_{[nt]}$. Section 3 contains the statements and the discussions of the main results, and Section 4 contains the proofs of them as consequences of (1). Section 5 contains auxiliary results for establishing the approximation (1). The proof of (1) constitutes Section 6.

NOTATION. The constants $\gamma_n$, $\beta_n$ and the functions $f(y)$ and $f_n(y) = \beta_n f(\beta_n y)$ are exclusively used in the sense they are used above. Similarly, $\phi(z)$ is used in the sense of (1) above, $b_n$ in the sense of (4), and $M_{h,\eta}(y)$ and $m_{h,\eta}(y)$ in the sense of (6). We let

$$\phi_\varepsilon(z) = \frac{1}{\varepsilon\sqrt{2\pi}} \exp\left\{ -\frac{z^2}{2\varepsilon^2} \right\}$$

[so that $\phi(z) = \phi_1(z)$]. $\widehat{K}$ stands for the Fourier transform of the measure $K$, $\mathbb{I}_A(\cdot)$ stands for the indicator function of the set $A$ and $\mathbb{R}$ stands for the real line. Convergence in $\mathbf{L}^2$ has the usual meaning of the convergence in mean-square. The notation $C$ stands for a generic constant that may take different values even at different places of the same proof.

## 2. Preliminaries.
Let $\{Z_\alpha(t), t \in \mathbb{R}\}$, $0 < \alpha \leq 2$ be an $\alpha$-stable Lévy motion. This means $Z_\alpha(t)$ has stationary independent increments having a strictly $\alpha$-stable distribution, that is, for $s < t$, $Z_\alpha(t) - Z_\alpha(s)$ has the characteristic function $\exp\{-(t-s)|u|^\alpha(1 + i\beta \, \mathrm{sign}(u) \tan(\frac{\pi\alpha}{2}))\}$, where $|\beta| \leq 1$ with $\beta = 0$ when $\alpha = 1$. (Note that this definition of strict $\alpha$-stability for the case $\alpha = 1$ differs from the usual one in that we take the shift parameter to be 0.) When $\alpha = 2$, $Z_\alpha(t)$ becomes the Brownian motion with variance 2.

A process $\{\Lambda_{\alpha,H}(t), t \geq 0\}$ is called a LFSM with Hurst parameter $H$, $0 < H < 1$, if it is given by

$$\Lambda_{\alpha,H}(t) = a \int_{-\infty}^{0} \{(t-u)^{H-1/\alpha} - (-u)^{H-1/\alpha}\} Z_\alpha(du)$$

$$+ a \int_{0}^{t} (t-u)^{H-1/\alpha} Z_\alpha(du),$$

where $Z_\alpha(t)$ is an $\alpha$-stable Lévy motion as above and $a$ is a nonzero constant. When $\alpha = 2$, the LFSM reduces to the FBM. See Samorodnitsky and Taqqu (1994) and Maejima (1989) for the details of LFSM.

We make the convention that in the case $H = 1/\alpha$, the LFSM $\{\Lambda_{\alpha,H}(t), t \geq 0\}$ is taken to be $\{Z_\alpha(t), t \geq 0\}$. It is important to note, however, that in this case the restriction $0 < H < 1$ is equivalent to that of $1 < \alpha \leq 2$.

Let $\{\zeta(t), t \geq 0\}$ be a real valued measurable process. Then a measurable process $\{L(t,x), t \geq 0, x \in \mathbb{R}\}$ is said to be a *local time* of $\{\zeta(t), t \geq 0\}$ if for each $t \geq 0$,

$$(2) \qquad \int_0^t \mathbb{I}_A(\zeta(s))\,ds = \int_{-\infty}^\infty \mathbb{I}_A(x)L(t,x)\,dx \qquad \text{for all Borel subset } A \text{ of } \mathbb{R}$$

with probability one. [Without loss of generality we take $\Lambda_{\alpha,H}(t)$ to be measurable.] For the symmetric LFSM, the existence of $L(t,x)$ is known, see Kôno and Maejima (1991). For the general LFSM, we have the following result. [It may be noted that when $0 < \alpha \leq 1$, the local time for $\{Z_\alpha(t)\}$ does not exist; the case $1 < \alpha \leq 2$ is covered by the next result.]

THEOREM 0.  *For a LFSM $\{\Lambda_{\alpha,H}(t), t \geq 0\}$ with $0 < \alpha \leq 2$ and $0 < H < 1$, there is a local time $L(t,x)$ such that for each $t$ and $x$,*

$$\lim_{\eta \downarrow 0} \frac{1}{\eta} \int_0^t \mathbb{I}_{[x, x+\eta)}(\Lambda_{\alpha,H}(s))\,ds = L(t,x) \qquad in \; \mathbf{L}^2.$$

*In addition, $L(t,x)$ has the representation $L(t,x) = \frac{1}{2\pi} \int_{-\infty}^\infty \int_0^t e^{iu(\Lambda_{\alpha,H}(s)-x)}\,ds\,du$.*

Next we recall the result on the weak convergence of the partial sum process $S_{[nt]}$. Recall that $\xi_j, -\infty < j < \infty$, is a sequence of i.i.d. r.v.s belonging to the domain of attraction of a strictly stable law with index $0 < \alpha \leq 2$. For the later purpose we mention that this, in particular, means for all $u$ in some neighborhood of 0,

$$(3) \qquad E[e^{iu\xi_1}] = \begin{cases} \exp\left\{-|u|^\alpha G(|u|)\left(1 + i\beta\,\mathrm{sign}(u)\tan\left(\dfrac{\pi\alpha}{2}\right)\right)\right\}, & \text{if } \alpha \neq 1, \\ \exp\{-|u|G(|u|)\}, & \text{if } \alpha = 1, \end{cases}$$

with $|\beta| \leq 1$, where $G(u)$ is slowly varying as $u \to 0$. In addition, if one lets

$$(4) \qquad b_n^{-1} = \inf\{u > 0 : u^\alpha G(u) = n^{-1}\},$$

then $b_n$ is of the form $n^{1/\alpha}h(n)$ for some slowly varying $h(n)$; in fact, $b_n^\alpha \backsim nG(b_n^{-1})$. [For the details of these facts, see, e.g., Bingham, Goldie and Teugels (1987), page 344.]

Now recall that $c_j, j = 0, 1, \ldots,$ is a sequence of real numbers such that $c_0 = 1$, and $X_k = \sum_{j=0}^\infty c_j \xi_{k-j}$. We impose the following (mutually exclusive) conditions:

(C.1)

(5) $$c_j = j^{H-1-1/\alpha} R(j) \qquad \text{with } H \neq 1/\alpha, 0 < H < 1,$$

where $R(j)$ is slowly varying at infinity, and

$$\sum_{j=0}^{\infty} c_j = 0 \qquad \text{when } H - 1/\alpha < 0.$$

(C.2)

$$\sum_{j=0}^{\infty} |c_j|^\tau < \infty \qquad \text{for some } \tau \text{ such that } 0 < \tau < \alpha, \ \tau \leq 1$$

and

$$\sum_{j=0}^{\infty} c_j \neq 0.$$

Throughout what follows we let, with $b_n$ as in (4),

$$\gamma_n = \begin{cases} n^{H-1/\alpha} R(n) b_n, & \text{if the condition (C.1) is satisfied,} \\ \left( \sum_{j=0}^{\infty} c_j \right) b_n, & \text{if the condition (C.2) is satisfied.} \end{cases}$$

Because $b_n = n^{1/\alpha} h(n)$, one has $n^{H-1/\alpha} R(n) b_n = n^H u(n)$ for a slowly varying $u(n)$.

The following result is taken from Kasahara and Maejima [(1988), Theorems 5.1–5.3], but see also the references given there, especially Astrauskas (1983).

PROPOSITION 1.   (i) *Assume that the condition* (C.1) *is satisfied. Then the finite-dimensional distributions of the process* $\gamma_n^{-1} S_{[nt]}$ *converge in distribution to those of the LFSM process* $\Lambda_{\alpha,H}(t)$, $H \neq 1/\alpha$.

(ii) *Assume that the condition* (C.2) *is satisfied. Then the finite-dimensional distributions of the process* $\gamma_n^{-1} S_{[nt]}$ *converge in distribution to those of the* $\alpha$-*stable Lévy motion* $Z_\alpha(t)$.

Note that the statement (ii) of this result holds for the entire range $0 < \alpha \leq 2$, but our interest is only on the range $1 < \alpha \leq 2$ because, as was mentioned earlier, the local time for $Z_\alpha(t)$ does not exist when $0 < \alpha \leq 1$. So in the situation of (5) with $H - 1/\alpha < 0$ but $\sum_{j=0}^{\infty} c_j \neq 0$ [which case was excluded in the statement (i)], this restricts $\alpha$ to either $1 < \alpha < \min\{2, \frac{1}{H}\}$ or $\alpha = 2$ when $\frac{1}{H} > 2$.

**3. Statements and discussions of the main results.** It is assumed that the constants $\beta_n, n \geq 1$, involved throughout below is such that $\beta_n \to \infty$ with $n^{-1}\beta_n \to 0$ as $n \to \infty$. For any function $h(y)$, we define

(6)  $M_{h,\eta}(y) = \sup\{h(u) : |u - y| \leq \eta\}$  and  $m_{h,\eta}(y) = \inf\{h(u) : |u - y| \leq \eta\}.$

Also, under the condition (C.2), we shall henceforth restrict to the situation $1 < \alpha \leq 2$ for the reason mentioned above, so that according to our convention the limit $Z_\alpha(t)$ in Proposition 1(ii) becomes $\Lambda_{\alpha,H}(t)$ with $H = 1/\alpha$.

THEOREM 2. *Assume that either one of the conditions* (C.1) *or* (C.2) *is satisfied. In addition, assume that* $\xi_1$ *satisfies Cramér's condition*

$$\limsup_{|u| \to \infty} |E[e^{iu\xi_1}]| < 1.$$

*Let,* $f(y)$ *be such that* $M_{|f|,\eta}(y)$ *and* $M_{f^2,\eta}(y)$ *are Lebesgue integrable for some* $\eta > 0$ *and*

(7)  $$\int_{-\infty}^{\infty} (M_{f,\delta}(y) - m_{f,\delta}(y))\,dy \to 0 \qquad as \ \delta \to 0.$$

*Then the finite-dimensional distributions of the process* $n^{-1}\beta_n \sum_{k=1}^{[nt]} f(\beta_n \times (\gamma_n^{-1}S_k + x))$ *[indexed by* $(t, x)$*] converge in distribution to those of* $(\int_{-\infty}^{\infty} f(y)\,dy)L(t, -x)$.

Note that $L(t, -x)$ will have the same distribution as that of $L(t, x)$ only when the $\alpha$-stable Lévy motion $Z_\alpha(t)$ involved in the definition of LFSM is symmetric around zero (which is always true in the case of FBM, the case $\alpha = 2$).

REMARK 1. There are alternative requirements on $f(y)$ that will imply those stated in Theorem 2. For example, one possibility is to assume that the set of discontinuity points of $f(y)$ is of Lebesgue measure zero, together with the Lebesgue integrability of $M_{|f|,\eta}(y)$. [It is clear that condition (7) is then implied by the dominated convergence theorem]. Also, as will be indicated later (see the Remark 7), it is possible to relax the Lebesgue integrability of $M_{f^2,\eta}(y)$ to that of local Lebesgue integrability. Thus the second possibility is to assume the local Riemann integrability of $f(y)$, together with the Lebesgue integrability of $M_{|f|,\eta}(y)$. In particular, Theorem 2 holds for the important situation in which $f(y) = I_{(c,d)}(y)$. (Here the limit will remain the same if the open interval $(c, d)$ is replaced by the closed interval $[c, d]$ or by a semi-open interval.)

REMARK 2. In the i.i.d. situation $S_k = \sum_{j=1}^k \xi_j$, with $\beta_n = \gamma_n = b_n = \sqrt{n}$, and when $f(y)$ is assumed to be Riemann integrable such that $|f(y)| \leq B(1 + |y|^{-1-c})$ for some $B > 0$ and $c > 0$, Borodin and Ibragimov [(1995),

Theorem 2.1, Chapter IV, page 143, and Theorem 2.2, Chapter IV, page 145] show that the conclusion of Theorem 2 holds without the Cramér's condition. When $f(y)$ is as above, we mention, without going into the details, the following partial extensions: (a) Theorem 2 extends to the situation where $\xi_1$ is nonlattice without satisfying the Cramér's condition but only under the restriction $b_n^{-1}\beta_n = O(1)$, where $b_n$ is as in (4). [In the important case $\beta_n = \gamma_n$, this will include the situation of condition (C.2), as well as that of condition (C.1) when $H - 1/\alpha < 0$, but unfortunately will exclude the case $H - 1/\alpha > 0$.] (b) Similarly, Theorem 2 extends to the situation where $\xi_1$ has a lattice distribution when $b_n^{-1}\beta_n \to 0$ [which in the case $\beta_n = \gamma_n$ will include the situation of condition (C.1) when $H - 1/\alpha < 0$, but not otherwise].

REMARK 3.   Regarding the results available in the direction of this paper, Akonom (1993) deals with the situation (of a Brownian motion limit) where $\sum_{j=0}^{\infty} j |c_j| < \infty$ with $\sum_{j=0}^{\infty} c_j \neq 0$, $E[|\xi_1|^p] < \infty$ for some $p > 2$, $\lim_{u \to \infty} u^\tau E[e^{iu\xi_1}] = 0$ for some $\tau > 0$ and the distribution of $\xi_1$ has a Lebesgue density. Then the main result obtained there implies the conclusion of Theorem 2 when $f(y) = I_{[a,b]}(y)$ and $\gamma_n = \sqrt{n}$. This situation is a special case of that of condition (C.2) for $\alpha = 2$, and Theorem 2 requires only that $\xi_1$ satisfies Cramér's condition and is in the domain of attraction of a normal distribution and that $\sum_{j=0}^{\infty} |c_j| < \infty$ with $\sum_{j=0}^{\infty} c_j \neq 0$. The underlying situation (with a Brownian motion limit) of Park and Phillips (1999, 2001) is the same as that of Akonom (1993) but the form of $f(y)$ is not restricted to $f(y) = I_{[a,b]}(y)$, however, the (Lipschitz type) conditions imposed there on $f(y)$ are unfortunately rather strong in addition to requiring further moment conditions on $\xi_1$.

Tyurin and Phillips (1999) consider the situation where $X_k$ is in the truncated form $X_k = \sum_{j=0}^{k} c_j \xi_{k-j}$, in addition to the restriction $H - 1/2 < 0$. In addition to the underlying assumptions of Akonom (1993) on the distribution of $\xi_1$ indicated above, it is further required that $E[|\xi_1|^q] < \infty$ for $q > \frac{1}{H}$. The limiting Gaussian process involved here will be different from a FBM, but will have similar properties. (It is easy to see that the results of the present paper hold for the truncated case $X_k = \sum_{j=0}^{k} c_j \xi_{k-j}$ also with the changes in the limiting forms taken into account.)

Now, an example given in Borodin and Ibragimov [(1995), Chapter IV, page 143] shows that the requirement (7) on $f(y)$ in Theorem 2 cannot be avoided entirely. The next result relaxes that requirement, but assumes conditions stronger than the Cramér's condition.

THEOREM 3.   *Assume that either one of the conditions* (C.1) *or* (C.2) *is satisfied.*

(i) *Suppose that, for some integer $n_0$, the $n_0$-fold convolution of the distribution of $\xi_1$ has a nonzero absolutely continuous component. Let $f(y)$ be Lebesgue integrable such that $\sup_{y \in \mathbb{R}} |f(y)| < \infty$. Then the conclusion of Theorem 2 holds.*

(ii) *Suppose that, for some integer $n_0 > 0$, $\int |\prod_{j=0}^{n_0-1} E[e^{iug(j)\xi_j}]| \, du < \infty$ where $g(j) = \sum_{i=0}^{j} c_i$ and $n^{-1}\beta_n \times \sum_{k=1}^{n_0-1} f(\beta_n(\gamma_n^{-1}S_k + x))$ converges in probability to 0. Assume further that $f(y)$ and $f^2(y)$ are Lebesgue integrable. Then the conclusion of Theorem 2 holds.*

Note that if $\int |E[e^{iu\xi_1}]|^p du < \infty$ for some $p > 0$, then the condition $\int |\prod_{j=0}^{n_0-1} E[e^{iug(j)\xi_j}]| \, du < \infty$ for some integer $n_0 > 0$ in statement (ii) above is satisfied.

Note that the requirements on $f(y)$ in the first statement are stronger than those in the second statement; consider, for instance, $f(y) \sim |y|^\tau$ as $|y| \to 0$ with $0 > \tau > -1/2$. Also, as in Theorem 2 (see Remark 1), the integrability of $f^2(y)$ in the second statement can be relaxed to that of local integrability.

Let us make some remarks regarding the requirement on $n^{-1}\beta_n \sum_{k=1}^{n_0-1} f(\beta_n \times (\gamma_n^{-1}S_k + x))$ in statement (ii). First, it is redundant when $n_0 = 1$. When $x = 0$ and $\beta_n = \gamma_n$, an important case in applications, the requirement is satisfied because the quantity reduces to $n^{-1}\beta_n \sum_{k=1}^{n_0-1} f(S_k)$, which clearly converges in probability to 0 in view of $n^{-1}\beta_n \to 0$.

Suppose that $x \neq 0$. Then the additional condition $\lim_{a \to \infty} \sup_{|y| \geq a} |f(y)| < \infty$ is sufficient because, for each $1 \leq k < n_0$, with probability tending to one $\beta_n(\gamma_n^{-1}S_k + x)$ will be supported in a neighborhood of $\pm\infty$. The same is the case when $x = 0$ and $\beta_n\gamma_n^{-1} \to \infty$.

In the remaining case $x = 0$ and $\beta_n\gamma_n^{-1} \to 0$, suppose, for instance, that $f(y) \sim |y|^\tau$ as $|y| \to 0$ with $0 > \tau > -1/2$. Then with probability tending to one, $n^{-1}\beta_n|f(\beta_n\gamma_n^{-1}S_k)| \leq Cn^{-1}\beta_n|\beta_n\gamma_n^{-1}|^\tau = C(n^{-1}\beta_n)^{1+\tau}(n\gamma_n^{-1})^\tau \to 0$.

REMARK 4. So far our results are for $f_n(y) = \beta_n f(\beta_n y)$ based on $f(y)$. It is possible to extend Theorems 2 and 3 to more general $f_n(y)$ that satisfy the following conditions: $\sup_n \int |f_n(y)| \, dy < \infty$, $\sup_n \frac{1}{n} \int_{-\infty}^{\infty} |f_n^2(y)| \, dy \to 0$, $\lim_{\kappa \to \infty} \sup_n \int_{\{|y| \geq \kappa\}} |f_n(y)| \, dy = 0$ and, letting

$$F_n(y) = \begin{cases} \int_0^y f_n(u) \, du, & \text{if } y \geq 0, \\ -\int_y^0 f_n(u) \, du, & \text{if } y < 0, \end{cases}$$

there is an $F(y)$ such that $F_n(y) \to F(y)$ at all *continuity points* of $F(y)$. These conditions are satisfied when $f_n(y) = \beta_n f(\beta_n y)$ with $F(y) = \int_0^\infty f(u) \, du$ if $y \geq 0$ and $F(y) = \int_{-\infty}^0 f(u) \, du$ if $y < 0$, and, in fact, we shall use the assumptions on $f(y)$ only in the form of the above conditions. Theorem 3(ii) extends as follows. Assume, for some integer $n_0 > 0$, $\int |E[e^{iu\xi_1}]|^{n_0} \, du < \infty$ and $n^{-1} \sum_{k=1}^{n_0-1} f_n(\gamma_n^{-1} S_k + x)$ converges in probability to 0. Then, if $f_n(y)$ satisfy the above conditions, $n^{-1} \sum_{k=1}^{[nt]} f_n(\gamma_n^{-1} S_k + x)$ converges in distribution to $\int_{-\infty}^{\infty} L(t, y - x) \, dF(y)$. Similarly, Theorem 3(i) extends, under the additional assumption $\sup_{\{n,y\}} \sigma_n^{-1} |f_n(y)| < \infty$ for some $\sigma_n \to \infty$ with $\frac{\sigma_n}{n} \to 0$.

The next result may be viewed as a discrete approximation to the local time of the LFSM, which we obtain as a by-product to Theorem 3(ii) (whose requirements are satisfied with $n_0 = 1$). Note that the approximation is in $\mathbf{L}^2$, in contrast to the distributional approximation in Theorems 2 and 3. Note also that in this case $\gamma_n = n^H$.

THEOREM 4. *Suppose that $f(y)$ and $f^2(y)$ are Lebesgue integrable. Then, for each $t$,*

$$\frac{\beta_n}{n} \sum_{k=1}^{[nt]} f\left(\beta_n\left(\Lambda_{\alpha,H}\left(\frac{k}{n}\right) - x\right)\right) \to \left(\int_{-\infty}^{\infty} f(y) \, dy\right) L(t, x) \qquad in \ \mathbf{L}^2.$$

The next result is a continuous analogue of the preceding result.

THEOREM 5. *Assume that $f(y)$ and $f^2(y)$ are Lebesgue integrable. Then*

$$\frac{1}{\kappa^{1-H}} \int_0^{\kappa t} f\left(\kappa^H\left(\Lambda_{\alpha,H}\left(\frac{s}{\kappa}\right) - x\right)\right) ds \to \left(\int_{-\infty}^{\infty} f(y) \, dy\right) L(t, x) \qquad in \ \mathbf{L}^2$$

*for each $t$ and $x$ as $\kappa \to \infty$.*

As noted in connection with Theorem 3(ii), the Lebesgue integrability of $f^2(y)$ in Theorems 4 and 5 can be relaxed to that of local Lebesgue integrability.

REMARK 5. Note that because the distribution of $n^H \Lambda_{\alpha,H}(k/n)$ is the same as that of $\Lambda_{\alpha,H}(k)$, it follows from Theorem 4 that $\frac{1}{n^{1-H}} \sum_{k=1}^{[nt]} f(\Lambda_{\alpha,H}(k) - xn^H)$ converges in distribution to $(\int_{-\infty}^{\infty} f(y) \, dy) L(t, x)$ as $n \to \infty$. Similarly, $\frac{1}{\kappa^{1-H}} \int_0^{\kappa t} f(\Lambda_{\alpha,H}(s) - x\kappa^H) \, ds$ converges in distribution to $(\int_{-\infty}^{\infty} f(y) \, dy) L(t, x)$ as $\kappa \to \infty$, in view of Theorem 5.

REMARK 6. One has $\frac{1}{\kappa^{1-H}} \int_0^{\kappa t} f(\kappa^H(\Lambda_{\alpha,H}(s/\kappa) - x)) \, ds = \kappa^H \int_0^{\kappa t} f(\kappa^H \times (\Lambda_{\alpha,H}(s) - x)) \, ds$. When $f(y) = \mathbb{I}_{[0,1)}(y)$ and $\kappa^H = 1/\eta$, the right-hand side reduces to $\frac{1}{\eta} \int_0^t \mathbb{I}_{[x,x+\eta)}(\Lambda_{\alpha,H}(s)) \, ds$.

**4. Proofs of the results.** The proposition stated next is just a formalization of the approximation (1), the proof of which is postponed to the next two sections because it requires many auxiliary results. In this section we derive Theorems 0 and 2–5 as consequences of it. Recall that we let

$$\phi_\varepsilon(z) = \frac{1}{\varepsilon\sqrt{2\pi}} \exp\left\{ -\frac{z^2}{2\varepsilon^2} \right\} \quad \text{and} \quad \phi(z) = \frac{1}{\sqrt{2\pi}} \exp\left\{ -\frac{z^2}{2} \right\}.$$

PROPOSITION 6. *Let* $f_n(y) = \beta_n f(\beta_n y)$. *Assume that the conditions of any one of Theorems 2 or 3 are satisfied. Let* $L_n(t,x) = \frac{1}{n} \sum_{k=q}^{[nt]} f_n(\gamma_n^{-1} S_k + x)$, *where* $q = 1$, *except under the second part of Theorem 3 in which case* $q = n_0$ *with* $n_0$ *as involved there, and*

$$L_{n,\varepsilon}(t,x) = \frac{1}{n} \sum_{k=1}^{[nt]} \int_{-\infty}^{\infty} f_n(\gamma_n^{-1} S_k + x + z\varepsilon)\phi(z)\,dz.$$

*Then*

$$\lim_{\varepsilon \to 0} \lim_{n \to \infty} \sup_{t,x} E[(L_n(t,x) - L_{n,\varepsilon}(t,x))^2] = 0.$$

We next deal with $L_{n,\varepsilon}(t,x)$ through the following series of steps.

LEMMA 7. *For each* $\varepsilon > 0$,

$$(8) \qquad \sup_{t,x} \left| L_{n,\varepsilon}(t,x) - \left( \int_{-\infty}^{\infty} f(z)\,dz \right) \frac{1}{n} \sum_{k=1}^{[nt]} \phi_\varepsilon(\gamma_n^{-1} S_k + x) \right|$$

*is bounded by a nonrandom quantity that converges to 0 as* $n \to \infty$.

PROOF. For convenience take $\varepsilon = 1$. Let $F_n(y) = \int_{-\infty}^{y} \beta_n f(\beta_n u)\,du$, and define $F(y)$ such that $F(y) = \int_{-\infty}^{\infty} f(u)\,du$ if $y \geq 0$ and $F(y) = 0$ if $y < 0$. Note that $F_n(y) \to F(y)$ at all continuity points of $F(y)$ and $F(b) - F(a) = 0$ if $0 \notin (a,b]$. Now $L_{n,\varepsilon}(t,x)$ takes the form $\int_{-\infty}^{\infty} (\frac{1}{n} \sum_{k=1}^{[nt]} \phi(y - \gamma_n^{-1} S_k - x))\,dF_n(y)$. The difference between this and

$$(9) \qquad \int_{\{|y| \leq \nu\}} \left( \frac{1}{n} \sum_{k=1}^{[nt]} \phi(y - \gamma_n^{-1} S_k - x) \right) dF_n(y)$$

is bounded in absolute value by $C \int_{\{|y| > \nu\}} d|F_n|(y) = C \int_{\{|u| > \beta_n \nu\}} |f(u)|\,du$. Define $y_{m,i}$, $i = -m, \ldots, 0, \ldots, m$ such that $y_{m,-m} = -\nu < y_{m,-m+1} < \cdots <$

$y_{m,m-1} < y_{m,m} = \nu$ and $\sup_i |y_{m,i} - y_{m,i-1}| \le \frac{2\nu}{m}$. Then the difference between (9) and

$$(10) \qquad \sum_{i=-m}^{m-1} \left( \frac{1}{n} \sum_{k=1}^{[nt]} \phi(y_{m,i} - \gamma_n^{-1} S_k - x) \right) \int_{y_{m,i}}^{y_{m,i+1}} dF_n(y)$$

is bounded in absolute value by $C \frac{\nu}{m} \int_{\{|y| \le \nu\}} d|F_n|(y) \le C \frac{\nu}{m}$. Further, the difference between (10) and $\sum_{i=-m}^{m-1} (\frac{1}{n} \sum_{k=1}^{[nt]} \phi(y_{m,i} - \gamma_n^{-1} S_k - x)) \int_{ym,i}^{y_{m,i+1}} dF(y)$ is bounded in absolute value by $C \sum_{i=-m}^{m-1} |\int_{y_{m,i}}^{y_{m,i+1}} d(F_n(y) - F(y))|$. Thus, it follows that (8) is bounded in absolute value by a constant multiple of $\int_{\{|u| > \beta_n \nu\}} |f(u)| \, du + \frac{\nu}{m} + \sum_{i=-m}^{m-1} |\int_{y_{m,i}}^{y_{m,i+1}} d(F_n(y) - F(y))|$. Denote this by $Q(\nu, m, n)$. It is clear that $\lim_{\nu \to \infty} \lim_{m \to \infty} \lim_{n \to \infty} Q(\nu, m, n) = 0$. $\quad \square$

LEMMA 8. *For each $\varepsilon > 0$, the finite-dimensional distributions of $\frac{1}{n} \times \sum_{k=1}^{[nt]} \phi_\varepsilon(\gamma_n^{-1} S_k + x)$ [indexed by $(t, x)$] converge in distribution to those of $\int_0^t \phi_\varepsilon(\Lambda_{\alpha, H}(s) + x) \, ds$.*

PROOF. For notational convenience, take $\varepsilon = 1$. Also, we take $x = 0$ for simplicity so that we consider the process $n^{-1} \sum_{k=1}^{[nt]} \phi(\gamma_n^{-1} S_k)$ indexed only by $t$. We now invoke Gikhman and Skorokhod [(1969), Theorem 1, page 485]. According to this result, letting $\zeta_n(s) = \phi(\gamma_n^{-1} S_{[ns]})$, we need to check that [the finite-dimensional distributions of $\zeta_n(s)$ being convergent to those of the limit $\zeta(s) = \phi(\Lambda_{\alpha, H}(s))$ by Proposition 1]

$$\sup_{s,n} E[|\zeta_n(s)|] < \infty, \qquad \lim_{\eta \to 0} \limsup_{n \to \infty} \sup_{|s_1 - s_2| \le \eta} E[|\zeta_n(s_1) - \zeta_n(s_2)|] = 0.$$

Because $\lim_{\kappa \to \infty} \sup_{|y| > \kappa} \phi(y) = 0$, it is enough to verify that, for every $\delta > 0$,

$$\lim_{\eta \to 0} \limsup_{n \to \infty} \sup_{|s_1 - s_2| \le \eta} P(|\gamma_n^{-1} S_{[ns_1]} - \gamma_n^{-1} S_{[ns_2]}| > \delta) = 0.$$

Note that $\sup_{|s_1 - s_2| \le \eta} P(|\gamma_n^{-1} S_{[ns_1]} - \gamma_n^{-1} S_{[ns_2]}| > \delta) = \sup_{0 < \tau \le \eta} P(|\gamma_n^{-1} \times S_{[n\tau]}| > \delta)$ by the stationarity of $\{X_k, k \ge 1\}$. (Here the stationarity is used only for convenience and can be avoided.) Now, Theorems 2.2 and 4.1 in Kasahara and Maejima (1988), give bounds for $\sup_{0 < \tau \le \eta} P(|\gamma_n^{-1} S_{[n\tau]}| > \delta)$, respectively, when $0 < \alpha < 2$ and when $\alpha = 2$ . Using the arguments similar to those used in Kasahara and Maejima [(1988), Section 5], these bounds converge to zero by first letting $n \to \infty$ and then $\eta \to 0$. This completes the proof. $\quad \square$

The next result is related to Theorem 4.

LEMMA 9. *For each $\varepsilon > 0$ and for all $t$ and $x$, $n^{-1} \sum_{k=1}^{[nt]} \phi_\varepsilon(\Lambda_{\alpha, H}(\frac{k}{n}) + x) \to \int_0^t \phi_\varepsilon(\Lambda_{\alpha, H}(s) + x) \, ds$ in $\mathbf{L}^2$.*

PROOF.    $t, x$ and $\varepsilon$ being fixed, take for notational simplicity $t = 1$, $x = 0$ and $\varepsilon = 1$. Write $n^{-1} \sum_{k=1}^{n} \phi(\Lambda_{\alpha,H}(\frac{k}{n}))$ in the form $\int_0^1 \phi(\Lambda_{\alpha,H}(\frac{[ns]}{n})) \, ds$. Then the proof is clear because $\Lambda_{\alpha,H}(\frac{[ns]}{n}) - \Lambda_{\alpha,H}(s)$ converges to 0 in probability [see Samorodnitsky and Taqqu 1994, Proposition 7.4.3] and because $\sup_z \phi(z) < \infty$. $\qquad \square$

The next result is a continuous analogue of Proposition 6 for the LFSM situation.

PROPOSITION 10.    *Suppose that $f(y)$ and $f^2(y)$ are Lebesgue integrable. Let $L_\kappa^*(t,x) = \frac{1}{\kappa^{1-H}} \int_0^{\kappa t} f(\kappa^H(\Lambda_{\alpha,H}(s/\kappa) - x)) \, ds$. Further, let $L_\varepsilon(t,x) = allow$ ($\int_{-\infty}^{\infty} f(y) \, dy) \int_0^t \phi_\varepsilon(\Lambda_{\alpha,H}(s) - x) \, ds$. Then*

$$\lim_{\varepsilon \to 0} \lim_{\kappa \to \infty} \sup_{x,t} E[|L_\kappa^*(t,x) - L_\varepsilon(t,x)|^2] \to 0.$$

The proof of this result is essentially contained in the proof of Proposition 6 for the situation of the second part (with $n_0 = 1$) of Theorem 3. The next proposition will be the only remaining fact required (apart from the proof of Proposition 6) to complete the proofs of Theorems 0 and 2–5.

PROPOSITION 11.    *For each $t$ and $x$, $\int_0^t \phi_\varepsilon(\Lambda_{\alpha,H}(s) - x) \, ds \to L(t,x)$ in $\mathbf{L}^2$ as $\varepsilon \downarrow 0$, where $L(t,x)$ is the local time of $\Lambda_{\alpha,H}(s)$, that is, satisfies the requirement (2). In addition, $L(t,x)$ has the representation of Theorem 0, that is, $L(t,x) = \frac{1}{2\pi} \int_{-\infty}^{\infty} \int_0^t e^{iu(\Lambda_{\alpha,H}(s) - x)} \, ds \, du$.*

PROOF.    Let $\Upsilon_\varepsilon = \int_0^t \phi_\varepsilon(\Lambda_{\alpha,H}(s) - x) \, ds$. We first show that $\Upsilon_\varepsilon$ is a Cauchy sequence in $\mathbf{L}^2$. For notational convenience, we restrict for the moment to $t = 1$ and suppress the occurrence of $x$. We need to show that $E[(\Upsilon_{\varepsilon_1} - \Upsilon_{\varepsilon_2})^2] \to 0$ as $(\varepsilon_1, \varepsilon_2) \downarrow 0$. Note that

$$e^{-v^2/2\varepsilon^2} = \frac{\varepsilon}{\sqrt{2\pi}} \int_{-\infty}^{\infty} e^{ivu} e^{-u^2 \varepsilon^2/2} \, du = \frac{\varepsilon}{\sqrt{2\pi}} \int_{-\infty}^{\infty} e^{-ivu} e^{-u^2 \varepsilon^2/2} \, du.$$

Hence, because $E[\Upsilon_{\varepsilon_1} \Upsilon_{\varepsilon_2}] = E[\int_0^1 \int_0^1 \phi_{\varepsilon_1}(\Lambda_{\alpha,H}(s_1)) \phi_{\varepsilon_2}(\Lambda_{\alpha,H}(s_2)) \, ds_1 \, ds_2]$,

$$\begin{aligned}
E[\Upsilon_{\varepsilon_1} \Upsilon_{\varepsilon_2}] &= \frac{1}{(2\pi)^2} \int_{-\infty}^{\infty} \int_{-\infty}^{\infty} E\left[ \int_0^1 \int_0^1 e^{iu\Lambda_{\alpha,H}(s_1) - iv\Lambda_{\alpha,H}(s_2)} \, ds_1 \, ds_2 \right] \\
&\qquad \times e^{-u^2 \varepsilon_1^2/2 - v^2 \varepsilon_2^2/2} \, du \, dv.
\end{aligned}$$
(11)

Here the order of integration is interchanged, which is permissible for each fixed $\varepsilon_1$ and $\varepsilon_2$, as can be seen using the fact $\int_0^1 \int_0^1 |E[e^{iu\Lambda_{\alpha,H}(s_1) - iv\Lambda_{\alpha,H}(s_2)}]| \, ds \, dt \leq$

1. Now if

$$(12) \quad \frac{1}{(2\pi)^2} \int_{-\infty}^{\infty} \int_{-\infty}^{\infty} E\left[\int_0^1 \int_0^1 e^{iu\Lambda_{\alpha,H}(s_1) - iv\Lambda_{\alpha,H}(s_2)} \, ds_1 \, ds_2\right] du \, dv < \infty,$$

then it is clear from (11) that $E[\Upsilon_{\varepsilon_1} \Upsilon_{\varepsilon_2}]$ converges to (12) as $(\varepsilon_1, \varepsilon_2) \to 0$ by the dominated convergence theorem. Note that

$$E\left[\int_0^1 \int_0^1 e^{iu\Lambda_{\alpha,H}(s_1) - iv\Lambda_{\alpha,H}(s_2)} \, ds_1 \, ds_2\right] \geq 0.$$

The same is also true for $E[\Upsilon_{\varepsilon_1}^2]$ and $E[\Upsilon_{\varepsilon_2}^2]$, so that $E[(\Upsilon_{\varepsilon_1} - \Upsilon_{\varepsilon_2})^2] \to 0$ as $(\varepsilon_1, \varepsilon_2) \downarrow 0$.

We now verify (12) in a stronger form that will be needed below. Note that, making now the occurrence of $x$ explicit, for the $L_\kappa^*(t,x)$ defined in Proposition 10, one can directly check that $\sup_\kappa \sup_x E[(L_\kappa^*(t,x))^2] < \infty$ (see Remark 8). Thus, in view of Proposition 10, it follows that $\limsup_{\varepsilon \to 0} \sup_x E[(\int_0^1 \phi_\varepsilon(\Lambda_{\alpha,H}(s)) \, ds)^2] < \infty$. Hence, in view of the identity (11) (with $\varepsilon_1 = \varepsilon_2$), it follows by Fatou's lemma that

$$
\begin{aligned}
&\sup_x \frac{1}{(2\pi)^2} \\
(13) \quad &\times \int_{-\infty}^{\infty} \int_{-\infty}^{\infty} \left[\int_0^1 \int_0^1 E[e^{iu(\Lambda_{\alpha,H}(s_1) - x) - iv(\Lambda_{\alpha,H}(s_2) - x)}] \, ds_1 \, ds_2\right] du \, dv \\
&< \infty.
\end{aligned}
$$

Actually we have shown that (11) converges to (12) because (11) is bounded above by (12).

Thus, we have established the Cauchy convergence of $\int_0^t \phi_\varepsilon(\Lambda_{\alpha,H}(s) - x) \, ds$ in $\mathbf{L}^2$, which entails the convergence in $\mathbf{L}^2$. Denote the limit by $L(t,x)$. In particular, (13) gives

$$(14) \quad \int_c^d \int_0^t \phi_\varepsilon(\Lambda_{\alpha,H}(s) - x) \, ds \, dx \to \int_c^d L(t,x) \, dx \qquad \text{in } \mathbf{L}^1$$

for every interval $[c,d]$. It remains to show that the limit $L(t,x)$ is indeed a local time, that is, (2) holds. It is sufficient to verify that $\int_0^t \mathbb{I}_{[c,d)}(\Lambda_{\alpha,H}(s)) \, ds = \int_c^d L(t,x) \, dx$ for every interval $[c,d)$. In view of Remark 6, Proposition 10, in particular, entails that

$$\frac{1}{\eta} \int_0^t \mathbb{I}_{[x,x+\eta)}(\Lambda_{\alpha,H}(s)) \, ds - \int_0^t \phi_\varepsilon(\Lambda_{\alpha,H}(s) - x) \, ds \to 0 \qquad \text{in } \mathbf{L}^2$$

as $\eta \to 0$ first and then $\varepsilon \to 0$, uniformly in $x$ (in the sense of Proposition 10). Hence,

$$\int_0^t \mathbb{I}_{[c,d)}(\Lambda_{\alpha,H}(s)) \, ds - \int_c^d \int_0^t \phi_\varepsilon(\Lambda_{\alpha,H}(s) - x) \, ds \, dx \to 0 \qquad \text{in } \mathbf{L}^2$$

as $\varepsilon \to 0$. In view of (14), this completes the proof.  $\square$

Now note that Propositions 10 and 11, together with Remark 6, proves Theorem 0. In the same way Theorems 2 and 3 follow from Propositions 6 and 11 and Lemmas 7 and 8. Thus, it only remains to establish Proposition 6.

**5. Auxiliary results for the proof of Proposition 6.** We first obtain some estimates on the behavior of the characteristic function of

$$(15) \qquad S_j^* = \gamma_j^{-1} \sum_{k=0}^{j-1} g(k) \xi_{j-k}.$$

(The reason for considering $S_j^*$, which is in the form of a finite-order sum of independent r.v.s, will become clear in the next section.) Here $g(k) = \sum_{i=0}^{k} c_i$ with $c_i$ as in the conditions (C.1) or (C.2). Only the distributional properties of $S_j^*$ will be required, so that for convenience we take $S_j^* = \gamma_j^{-1} \sum_{k=0}^{j-1} g(k) \xi_k$. Under the condition (C.1), we have $c_j = j^{H-1-1/\alpha} R(j)$ for a slowly varying $R(j)$ (with $\sum_{i=0}^{\infty} c_i = 0$ if $H - 1/\alpha < 0$). For simplicity, we shall restrict to the situation $c_j \sim C j^{H-1-1/\alpha}$, $j \to \infty$, for a suitable constant $C$. Then note that $g(k) = \sum_{i=0}^{k} c_i \sim C' k^{H-1/\alpha}$, $k \to \infty$, for some constant $C'$, so that for convenience $S_j^*$ can be taken to be of the form [with $b_j$ as in (4)]

$$(16) \qquad S_j^* = \sum_{k=1}^{j} \left(\frac{k}{j}\right)^{H-1/\alpha} \frac{\xi_k}{b_j}.$$

Note that in view of (3), there is a $\delta > 0$ such that, letting $\psi(v) = E[e^{iv\xi_1}]$,

$$(17) \qquad |\psi(v)| \le e^{-|v|^\alpha G(|v|)} \qquad \text{for all } |v| \le \delta,$$

where $G(w)$ is slowly varying at 0.

In the rest of this section it is assumed, without further mentioning, that either one of *the conditions* (C.1) *or* (C.2) *is satisfied.* Further, the constant $b_j$ is as in (4).

LEMMA 12. *Let $\widehat{H}_j(u)$ be the characteristic function of $S_j^*$ defined in* (15). *Then there are constants $A > 0, \lambda > 0,\ d > 0$ and $0 < c < \alpha$ such that $|\widehat{H}_j(u)| \le A e^{-d|u|^c}$ for all $|u| \le \lambda b_j$ and $j \ge 1$, where $b_j$ is as in* (4).

PROOF. It is enough to prove the result for all sufficiently large $j$, because for any $j_0 > 1$,

$$|\widehat{H}_j(u)| \le 1 = e^{d|\lambda b_{j_0}|^c} e^{-d|\lambda b_{j_0}|^c} \le e^{d|\lambda b_{j_0}|^c} e^{-d|u|^c}$$

$$\text{for } 1 \le j \le j_0 \text{ and } |u| \le \lambda b_j.$$

In the same way, it is enough to consider $u$ such that $|u| \geq C$ for some $C > 0$. We have $|\widehat{H}_j(u)| = \prod_{k=1}^j |\psi(\gamma_j^{-1} g(k)u)|$. Consider the situation of condition (C.1). Then [taking into account the simplification (16)], $|\widehat{H}_j(u)| = \prod_{k=1}^j |\psi((\frac{k}{j})^{H-1/\alpha} \frac{u}{b_j})|$. Suppose first that $H - 1/\alpha > 0$. Then $|(\frac{k}{j})^{H-1/\alpha} \frac{u}{b_j}| \leq \delta$ for all $|u| \leq \delta b_j$. Thus, in view of (17),

$$|\widehat{H}_j(u)| \leq \exp\left\{ -\sum_{k=1}^j \left|\left(\frac{k}{j}\right)^{H-1/\alpha} \frac{u}{b_j}\right|^\alpha G\left(\left(\frac{k}{j}\right)^{H-1/\alpha} \frac{|u|}{b_j}\right)\right\}$$

for all $|u| \leq \delta b_j$.

Recall that $b_j^\alpha \backsim jG(b_j^{-1})$ [the $G(u)$ in (4) and that in (17) being the same]. One can assume for convenience that $b_j^\alpha = jG(b_j^{-1})$. Then the sum in the preceeding exponent becomes

$$(18) \qquad |u|^\alpha \frac{1}{j} \sum_{k=1}^j \left|\left(\frac{k}{j}\right)^{H-1/\alpha}\right|^\alpha \frac{1}{G(b_j^{-1})} G\left(\left(\frac{k}{j}\right)^{H-1/\alpha} \frac{u}{b_j}\right).$$

Note that with $G(w)$ being continuous, it is bounded on compacts. Hence, by Potter's inequality [see Bingham, Goldie and Teugels (1987), statement (ii) of Theorem 1.5.6, page 25], for every $\delta > 0$ there is a $B > 0$, such that $|G(x)/G(y)| \leq B \max\{(x/y)^\delta, (x/y)^{-\delta}\}$ for all $x > 0$, $y > 0$. Hence, it follows easily that for a given $\epsilon > 0$, there are $C > 0$ and $C_1 > 0$ such that $\inf_{k \geq [j\epsilon]} \frac{1}{G(b_j^{-1})} G((\frac{k}{j})^{H-1/\alpha} \frac{|u|}{b_j}) \geq C_1 |u|^{-\delta}$ for all sufficiently large $j$ and for all $|u| \geq C$. Thus, (18) is bounded below, for some $C' > 0$, by $C'|u|^{\alpha-\delta}$ for all sufficiently large $j$ and for all $|u| \geq C$. This proves the result under condition (C.1) when $H - 1/\alpha > 0$.

When $H - 1/\alpha < 0$, note that whenever $k \geq [j\epsilon]$ for a given $\epsilon$, one has, for some $C > 0$, $(\frac{k}{j})^{H-1/\alpha} \leq C$ for all sufficiently large $j$. Hence, essentially the same arguments used above apply for this case also.

It remains to consider the case under condition (C.2). By assumption $\sum_{i=0}^\infty c_i \neq 0$ and $\sum_{i=0}^\infty |c_i| < \infty$, so that for some positive constants $\vartheta$ and $\kappa$ one has $\vartheta \leq |\sum_{i=0}^k c_i| \leq \kappa$ for all sufficiently large $k$. Further, in this case $\gamma_j^{-1} = b_j^{-1}(\sum_{i=0}^\infty c_i)^{-1}$. Hence, the proof of this case is contained in the above arguments. This completes the proof of the lemma. $\square$

LEMMA 13. *Let $\widehat{H}_j(u)$ be the characteristic function of $S_j^*$ defined in* (15). *Assume that $\xi_1$ satisfies the Cramér's condition of Theorem* 2. *Then for any $d > 0$, there is a $B > 0$ and a $0 < \rho < 1$ such that $\sup_{|u| \geq db_j} |\widehat{H}_j(u)| \leq B\rho^j$ for all $j \geq 1$, where $b_j$ is as in* (4).

PROOF. As in Lemma 12, it is enough to prove the statement for all sufficiently large $j$. Now recall that the Cramér's condition is actually equivalent to $\sup_{|z| \geq a} |\psi(z)| < 1$ for all $a > 0$. First consider the situation under (C.1). When $H - 1/\alpha < 0$ [and taking into account the simplification (16)], we have $(\frac{k}{j})^{H-1/\alpha} \geq 1$ for all $1 \leq k \leq j$. Thus,

$$\sup_{|u| \geq db_j} |\hat{H}_j(u)| = \sup_{|u| \geq db_j} \prod_{k=1}^{j} \left| \psi\left( \left(\frac{k}{j}\right)^{H-1/\alpha} \frac{u}{b_j} \right) \right| \leq \left( \sup_{|z| \geq d} |\psi(z)| \right)^j.$$

Now consider $H - 1/\alpha > 0$. Given $\epsilon > 0$, there is a $C > 0$ and a $j_0$ such that $(\frac{k}{j})^{H-1/\alpha} \geq C$ for all $k \geq [j\epsilon]$ and $j \geq j_0$. Hence, as before, $\sup_{|u| \geq db_j} |\hat{H}_j(u)| \leq (\sup_{|z| \geq dC} |\psi(z)|)^{j-[j\epsilon]}$, $j \geq j_0$. Because for some $\mu > 0$, $j - [j\epsilon] \geq \mu j$ for all sufficiently large $j$, this proves the result under (C.1). The proof under (C.2) uses the same arguments.   □

To proceed further we need the following result, contained in Bhattacharya and Ranga Rao [(1976), proof of Lemma 11.1, page 93].

LEMMA 14.   *Let $\mu$ and $\nu$ be finite measures on $R^k$. Let $\eta$ be a positive number and $K_\eta$ a probability measure on $R^k$ satisfying $K_\eta(\{x : |x| \leq \eta\}) = 1$. Let $h$ be a real valued Borel measurable function on $R^k$ such that $|M_{h,\eta}(x)|$ and $|m_{h,\eta}(x)|$ are integrable with respect to $\mu$ and $\nu$ [where $M_{h,\eta}(x)$ and $m_{h,\eta}(x)$ are as defined in (6)]. Then*

$$\int f\, d(\mu - \nu) \begin{cases} \geq \int m_{f,\eta}\, d(\mu - \nu) * K_\eta - \int (f - m_{f,2\eta})\, d\nu, \\ \leq \int M_{f,\eta}\, d(\mu - \nu) * K_\eta + \int (M_{f,2\eta} - f)\, d\nu. \end{cases}$$

The $K_\eta$ above will be chosen such that its characteristic function $\hat{K}_\eta(u)$ satisfies

$$(19) \qquad\qquad |\hat{K}_\eta(u)| \leq C \exp\{-(\eta|u|)^{1/2}\}$$

for all real $u$, where $C$ is a constant (independent of $\eta$). This is possible in view of Bhattacharya and Ranga Rao [(1976), Corollary 10.4, page 88], where $K_\eta$ is used extensively as a smoothing device. Next we give some inequalities that will be used in a rather crucial manner in the next section. For this reason we need to state the bounds explicitly.

LEMMA 15.   *Assume that $\xi_1$ satisfies the Cramér's condition. Let $K_\eta$ be the smoothing probability measure of Lemma 14 satisfying (18). Then, for some $0 < \rho < 1$,*

$$\int |\hat{H}_j(u)| |\hat{K}_\eta(u)|\, du \leq C(1 + \eta^{-1}\rho^j)$$

and

$$\int |e^{-u^2\varepsilon^2/2} - 1||\widehat{H}_j(u)||\widehat{K}_\eta(u)|\,du \le C(\varepsilon^2 + \eta^{-1}\rho^j).$$

PROOF. Let $\lambda > 0, \alpha > c > 0, d > 0$ be positive constants involved in Lemma 12 and let $0 < \rho < 1$ be as in Lemma 13. Then, recalling the inequality (19) for $|\widehat{K}_\eta(u)|$,

$$\int |\widehat{H}_j(u)||\widehat{K}_\eta(u)|\,du \le \int_{\{|u| \le \lambda b_j\}} |\widehat{H}_j(u)|\,du + \int_{\{|u| > \lambda b_j\}} |\widehat{H}_j(u)||\widehat{K}_\eta|\,du$$

$$\le C \int_{\{|u| \le \lambda b_j\}} e^{-d|u|^c}\,du + C\rho^j \int_{\{|u| > \lambda b_j\}} e^{-(\eta|u|)^{1/2}}\,du,$$

which proves the first part. The second part also follows in the same way using the inequality $|e^{-u^2\varepsilon^2/2} - 1| \le u^2\varepsilon^2/2$ over the range $\{|u| \le \lambda b_j\}$ and using the inequality $|e^{-u^2\varepsilon^2/2} - 1| \le 2$ over the range $\{|u| > \lambda b_j\}$. □

LEMMA 16. *Suppose that, for some integer $n_0$, the $n_0$-fold convolution of the distribution of $\xi_1$ has a nonzero absolutely continuous component. Let $H_j$ be the probability distribution $S_j^*$ defined in (15). Then there are measures $H_j^*$ and $H_j^{**}$ such that $H_j = H_j^* + H_j^{**}$, satisfying the following properties*:

(i) *For every $r > 0$, there is a constant $L_r$ such that $H_j^{**}(R) \le L_r j^{-r}$ for all $j \ge n_0$.*

(ii) *Let $\widehat{H}_j^*$ be the Fourier transform of $H_j^*$. Then there are constants $A > 0$, $\lambda > 0$, $d > 0$ and $0 < c < \alpha$ such that for every $j \ge n_0$ and for every $r > 0$, $|\widehat{H}_j^*| \le A e^{-d|u|^c} + L_r j^{-r}$ for all $|u| \le \lambda b_j$ and for a suitable constant $L_r$.*

(iii) *There are constants $B > 0$ and $0 < \rho < 1$ and an integrable function $g(v) \ge 0$ such that $\sup_{|u| \ge db_j} |\widehat{H}_j^*(u)| \le Bg(b_j^{-1}u)\rho^j$ for all $j \ge 2n_0$.*

PROOF. The proof is similar to the i.i.d. case $j^{-1/2}\sum_{i=1}^j \xi_i$ [see Prohorov (1952) and Le Cam (1960), pages 68–72]. We briefly recall the arguments. For simplicity, we consider the case $n_0 = 1$. First, letting $F$ for the distribution of $\xi_1$, one writes $F = \theta G_1 + (1-\theta)G_2$, where $0 < \theta < 1$ and $G_1$ and $G_2$ are probability measures such that $G_1$ is absolutely continuous with respect to the Lebesgue measure with a Lebesgue density uniformly bounded by a constant. (In particular, the density of $G_1$ is square integrable and, hence, $|\widehat{G}_1|^2$ is also integrable.) This is possible because $F$ is not entirely disjoint from the Lebesgue measure. Then the characteristic function $\psi_F(u)$ of $F$ has the corresponding decomposition $\psi_F(u) = \theta\psi_{G_1}(u) + (1-\theta)\psi_{G_2}(u)$. Hence, $\widehat{H}_j(u) = \prod_{k=1}^j (\theta\psi_{G_1}(\gamma_j^{-1}g(k)u) + (1-\theta)\psi_{G_2}(\gamma_j^{-1}g(k)u))$ which can be written in the form $\widehat{H}_j(u) = \sum_{l=0}^j \sum_{*l} \theta^l(1-\theta)^{j-l} \prod_{k=1}^j \psi_{*l}(\gamma_j^{-1}g(k)u)$, where $\psi_{*l}$ is either $\psi_{G_1}$ or $\psi_{G_2}$ such that $\psi_{G_1}$ occurs in the product $\prod_{k=1}^{j} \psi_{*l}(\gamma_j^{-1} g(k)u)$ exactly $l$ times, and the sum $\sum_{*l}$ is with respect to all such products for a given $l$. Now let $I_j^c = \{l : |l - j\theta| > j^{3/4}\}$. Then $H_j^{**}$ is the measure corresponding to the Fourier transform $\sum_{l \in I_j^c} \sum_{*l} \theta^l (1-\theta)^{j-l} \prod_{k=1}^{j} \psi_{*l}(\gamma_j^{-1} g(k)u)$. It is clear that $H_j^{**}(R)$ is bounded in absolute value by $\sum_{l \in I_j^c} \binom{j}{l} \theta^l (1-\theta)^{j-l}$, which is known to have the bound stated in statement (i). For statement (ii) use the fact that $|\hat{H}_j^*(u)| \le |\hat{H}_j(u)| + H_j^{**}(R)$ and then use the bound in statement (i) for $H_j^{**}(R)$, together with the bound in Lemma 12 for $|\hat{H}_j(u)|$.

Proof of statement (iii) uses essentially the same arguments of the i.i.d. case given in Prohorov (1952) or LeCam (1960) mentioned above, the essential facts being the Cramér's condition for $G_1$ and the function $g(v)$ in statement (iii) taken to be $|\hat{G}_1(av)|^2$ for a suitable constant $a$, which is integrable as indicated earlier. The necessary modifications needed for the present case being essentially the same as those used in the proofs of Lemmas 12 and 13, the proof is concluded.  □

Lemma 17.  *Under the situation of the preceeding lemma, for some $0 < \rho < 1$,*

$$\int |\hat{H}_j^*(u)| \, du \le C(1 + \rho^j)$$

*and*

$$\int |e^{-u^2 \varepsilon^2/2} - 1| |\hat{H}_j^*(u)| \, du \le C(\varepsilon^2 + \rho^j)$$

*for all $j \ge 2n_0$.*

Proof.  Regarding the first statement, according to statements (ii) and (iii) of Lemma 16, there is a $0 < \rho < 1$ such that for any $r > 0$, there is an $L_r$ satisfying

$$\int |\hat{H}_j^*(u)| \, du = \int_{\{|u| \le \lambda b_j\}} |\hat{H}_j^*(u)| \, du + \int_{\{|u| \ge \lambda b_j\}} |\hat{H}_j^*(u)| \, du$$
$$\le \int_{\{|u| \le \lambda b_j\}} (Ae^{-d|u|^c} + L_r j^{-r}) \, du + B\rho^j \int g(u) \, du$$

for all $j \ge 2n_0$. Here the integrable function $g(u) \ge 0$ is the one occurring in Lemma 16(iii). Similarly, the second statement follows.  □

The next result deals with the situation of the second statement of Theorem 3.

Lemma 18.  *Suppose that $\int |E[e^{iu\xi_1}]|^{n_0} \, du < \infty$ for some integer $n_0 > 0$. Then for all $j \ge n_0$, the conclusions of Lemma 17 hold for $\hat{H}_j(u)$ itself.*

PROOF. Because Lemma 16(iii) holds for $\widehat{H}_j(u)$ itself for all $j \geq n_0$ with $g(v) = |E[e^{iv\xi_1}]|^{n_0}$, the proof follows as in the previous lemma using Lemma 12.
$\qquad\square$

**6. Proof of Proposition 6.** We first prove the result under the conditions of Theorem 2 (assumed to hold henceforth without further mentioning), which will then essentially contain the proof under the conditions of Theorem 3. Throughout below the constants $\nu_n$ are such that

$$n^{-1}\nu_n \to 0 \quad \text{and} \quad \nu_n^{-1}\beta_n \to 0.$$

Such a choice of $\nu_n$ is possible because $n^{-1}\beta_n \to 0$. Further, $\rho$ will always stand for a constant with $0 < \rho < 1$. For the function $f_a(y) = f(ay)$, $a > 0$, the identities

$$
\begin{aligned}
(20) \qquad & \int M_{f_a, \eta/a}(y)\,dy = \frac{1}{a}\int M_{f,\eta}(y)\,dy, \\
& \int m_{f_a, \eta/a}(y)\,dy = \frac{1}{a}\int m_{f,\eta}(y)\,dy,
\end{aligned}
$$

which follow in view of $M_{f_a, \eta/a}(y) = M_{f,\eta}(ay)$ and $m_{f_a, \eta/a}(y) = m_{f,\eta}(ay)$, will be invoked repeatedly. In particular, because $\int M_{|f|,\eta}(y)\,dy < \infty$ for some $\eta > 0$,

$$(21) \qquad \limsup_{n\to\infty} \int M_{|f_n|,\nu_n^{-1}}(y)\,dy = \limsup_{n\to\infty} \int M_{|f|,\beta_n\nu_n^{-1}}(y)\,dy \leq C.$$

We begin with the following result, where and elsewhere recall that $f_n(y) = \beta_n f(\beta_n y)$.

PROPOSITION 19.

$$\sup_x E[|f_n(\gamma_n^{-1}S_k + x)|] \leq C(\gamma_n\gamma_k^{-1} + \nu_n\rho^k)\int M_{|f|,\beta_n\nu_n^{-1}}(y)\,dy.$$

PROOF. We have, with $g_n(j) = \gamma_n^{-1}\sum_{k=0}^{j} c_k$ [recall that $\gamma_n = n^{H-1/\alpha}R(n)b_n$],

$$\gamma_n^{-1}S_k = \sum_{j=-\infty}^{0}(g_n(k-j) - g_n(1-j))\xi_j + \sum_{j=1}^{k}g_n(k-j)\xi_j = S_{nk}^{**} + \gamma_n^{-1}\gamma_k S_k^*,$$

where $S_k^* = \sum_{j=1}^{k} g_k(k-j)\xi_j$. Note that $S_k^*$ and $S_{nk}^{**}$ are independent. Hence,

$$E[|f_n(\gamma_n^{-1}S_k + x)|] = E\left[\int |f_n(\gamma_n^{-1}\gamma_k y)|\,dH_k(y - \gamma_k^{-1}\gamma_n(x + S_{nk}^{**}))\right]$$

$$= \int |f_{nk}|\,dH_{*k},$$

where $H_k$ is the distribution of $S_k^*$, $H_{*k}(y) = H_k(y - \gamma_k^{-1}\gamma_n(x + S_{n,k}^{**}))$ and $f_{nk}(y) = f_n(\gamma_n^{-1}\gamma_k y)$. Noting that $S_k^*$ is in the form of (15), Lemma 14 (with the measure $v$ taken to be identically zero) gives $\int |f_{nk}| \, dH_{*k} \leq \int M_{|f_{nk}|, \nu_n^{-1}\gamma_n\gamma_k^{-1}} \times d(H_{*k} * K_{\nu_n^{-1}\gamma_n\gamma_k^{-1}})$. The right-hand side is bounded by $(\int M_{|f_{nk}|, \nu_n^{-1}\gamma_n\gamma_k^{-1}}(y) \, dy) \times (\int |\widehat{H}_k(u)| |\widehat{K}_{\nu_n^{-1}\gamma_n\gamma_k^{-1}}(u)| \, du)$ because for a probability measure $P$ with $\int |\widehat{P}(u)| \, du < \infty$, its density $p(y)$ satisfies

$$(22) \qquad \sup_y p(y) \leq \int |\widehat{P}(u)| \, du.$$

($H_{*k} * K_\eta$ has a density function.) Now $\int |\widehat{H}_k(u) \widehat{K}_{\nu_n^{-1}\gamma_n\gamma_k^{-1}}(u)| \, du \leq C(1 + \rho^k \nu_n \gamma_n^{-1}\gamma_k)$ according to Lemma 15. Also, $\int M_{|f_{nk}|, \nu_n^{-1}\gamma_n\gamma_k^{-1}}(y) \, dy = \gamma_n\gamma_k^{-1} \times \int M_{|f|, \beta_n\nu_n^{-1}}(y) \, dy$ in view of (20). Hence, the proof follows.    $\square$

Note that

$$(23) \qquad \sup_{\varepsilon, x} E\left[\int_{-\infty}^\infty |f_n(\gamma_n^{-1}S_k + x + z\varepsilon)|\phi(z) \, dz\right] \leq \sup_x E[|f_n(\gamma_n^{-1}S_k + x)|].$$

Further, $E[|f_n^2(\gamma_n^{-1}S_k + x)|] = \beta_n^2 E[|f^2(\beta_n(\gamma_n^{-1}S_k + x))|]$ and, hence, by Proposition 19,

$$\sup_x E[|f_n^2(\gamma_n^{-1}S_k + x)|] \leq C\beta_n(\gamma_n\gamma_k^{-1} + \nu_n\rho^k) \int M_{f^2, \beta_n\nu_n^{-1}}(y) \, dy.$$

Now note that $n^{-1}\sum_{k=1}^n (\gamma_n\gamma_k^{-1} + \nu_n\rho^k) \leq C$ because $\gamma_k = k^H u(k)$ with $0 < H < 1$ and $u(k)$ slowly varying. Thus, by (20) and because $n^{-1}\beta_n \to 0$,

$$(24) \qquad \left.\begin{array}{l} \sup_x \dfrac{1}{n^2} E\left[\displaystyle\sum_{k=1}^n |f_n^2(\gamma_n^{-1}S_k + x)|\right] \\[12pt] \sup_{\varepsilon, x} \dfrac{1}{n^2} E\left[\displaystyle\sum_{k=1}^n \int_{-\infty}^\infty |f_n^2(\gamma_n^{-1}S_k + x + z\varepsilon)|\phi(z) \, dz\right] \end{array}\right\} \to 0.$$

REMARK 7. If we let $f^{(\tau)}(y) = f(y)\mathbb{I}(|y| > \tau)$, then $\sup_x \beta_n E[|f^{(\tau)}(\beta_n \times (\gamma_n^{-1}S_k + x))|]$ is bounded by $C(\gamma_n\gamma_k^{-1} + \nu_n\rho^k) \int M_{|f^{(\tau)}|, \beta_n\nu_n^{-1}}(y) \, dy$, in view of Proposition 19, where $\lim_{\tau\to\infty} \lim\sup_{n\to\infty} \int M_{|f^{(\tau)}|, \beta_n\nu_n^{-1}}(y) \, dy = 0$ if $\int M_{|f|, \eta}(y) \, dy < \infty$ for some $\eta > 0$. Then without loss of generality, $f(y)$ in Theorem 2 can be taken to have a compact support; in particular, the integrability of $M_{f^2, \eta}(y)$ can be relaxed to local integrability.

The next lemma will ease the computations to be carried out further on.

Lemma 20. *Let $Z_k$, $k \geq 1$, be i.i.d. standard Gaussian random variables independent of $\{\xi_j : -\infty < j < \infty\}$. Then, for each $\varepsilon > 0$,*

$$\sup_{\varepsilon, t, x} E\left[\left\{\frac{1}{n}\sum_{k=1}^{[nt]}\left(f_n(\gamma_n^{-1}S_k + x + Z_k\varepsilon)\right.\right.\right.$$
$$\left.\left.\left. - \int_{-\infty}^{\infty} f_n(\gamma_n^{-1}S_k + x + z\varepsilon)\phi(z)\,dz\right)\right\}^2\right] \to 0.$$

Proof. Define the $\sigma$-fields $F_j = \sigma(S_1, \ldots, S_n, Z_1, \ldots, Z_j)$. Because $\{Z_k\}$ and $\{S_k\}$ are independent, $E[f_n(\gamma_n^{-1}S_k + x + Z_k\varepsilon)|F_{k-1}] = \int_{-\infty}^{\infty} f_n(\gamma_n^{-1}S_k + x + z\varepsilon)\phi(z)\,dz$. The differences $\frac{1}{n}\{f_n(\gamma_n^{-1}S_k + x + Z_k\varepsilon) - E[f_n(\gamma_n^{-1}S_k + x + Z_k\varepsilon)|F_{k-1}]\}$ form martingale differences with respect to $\{F_k, k \geq 1\}$. Hence, the expectation in the statement of the lemma is bounded by $E[\frac{1}{n^2}\sum_{k=1}^{n}\int_{-\infty}^{\infty} f_n^2(\gamma_n^{-1} \times S_k + x + z\varepsilon)\phi(z)\,dz]$. The proof follows by (24). □

Final arguments of the Proof of Proposition 6. For notational convenience, and because all the bounds derived below will be independent of $x$, we restrict to the case $x = 0$. In the same way, we restrict to $t = 1$. In view of the preceding Lemma 20, it is enough to show that $E[\{\frac{1}{n}\sum_{k=1}^{[nt]}(f_n(\gamma_n^{-1}S_k) - f_n(\gamma_n^{-1}S_k + Z_k\varepsilon))\}^2] \to 0$ by letting first $n \to \infty$ and then $\varepsilon \to 0$. Then, in view of (24), it is enough to show that

$$(25) \quad \frac{1}{n^2}\sum_{j=1}^{n}\sum_{k=j+1}^{n} E[(f_n(\gamma_n^{-1}S_j) - f_n(\gamma_n^{-1}S_j + Z_j\varepsilon))$$
$$\times (f_n(\gamma_n^{-1}S_k) - f_n(\gamma_n^{-1}S_k + Z_k\varepsilon))]$$

converges to 0 by first letting $n \to \infty$ and then $\varepsilon \to 0$. Recall that for $k > j$,

$$\gamma_n^{-1}S_k = \sum_{l=-\infty}^{0}(g_n(k-l) - g_n(1-l))\xi_l + \sum_{l=1}^{j} g_n(k-l)\xi_l + \sum_{l=j+1}^{k} g_n(k-l)\xi_l$$

with $g_n(j) = \gamma_n^{-1}\sum_{k=0}^{j} c_k$ and $\gamma_n = n^{H-1/\alpha}R(n)b_n$. This is the same as

$$\gamma_n^{-1}S_k = S_{nkj}^{**} + \gamma_n^{-1}\gamma_{k-j}S_{k-j}^{*},$$

where $S_{k-j}^{*} = \sum_{q=0}^{k-j-1} g_{k-j}(q)\xi_{k-q}$. Because $S_{k-j}^{*}$ and $S_{nkj}^{**}$ are independent,

$$E[f_n(\gamma_n^{-1}S_j)f_n(\gamma_n^{-1}S_k + Z_k\varepsilon)]$$
$$= E\left[f_n(\gamma_n^{-1}S_j)\int f_n(\gamma_n^{-1}\gamma_{k-j}y)\,d(H_{k-j} * G_{\varepsilon nkj})(y - \gamma_{k-j}^{-1}\gamma_n S_{nkj}^{**})\right],$$

where $H_{k-j}$ is the distribution of $S^*_{k-j}$, $G_\sigma$ is Gaussian with mean 0 and variance $\sigma^2$ and

$$\varepsilon_{nkj} = \varepsilon \gamma^{-1}_{k-j} \gamma_n.$$

Using similar identities, it follows that (25) can be written in the form

$$(26) \quad \frac{1}{n^2} \sum_{j=1}^{n} E\Bigg[ (f_n(\gamma^{-1}_n S_j) - f_n(\gamma^{-1}_n S_j + Z_j \varepsilon)) \\ \times \sum_{k=j+1}^{n} \int f_{nkj} \, d(H^{\#}_{k-j} - H^{\#}_{k-j} * G_{\varepsilon_{nkj}}) \Bigg],$$

where we have set $f_{nkj}(y) = f_n(\gamma^{-1}_n \gamma_{k-j} y)$ and $H^{\#}_{k-j}(y) = H_{k-j}(y - \gamma^{-1}_{k-j} \times \gamma_n S^{**}_{nkj})$. We first show that, for every $\delta > 0$, the quantity

$$(27) \quad \frac{1}{n^2} \sum_{j=[n\delta]}^{n} E\Bigg[ (f_n(\gamma^{-1}_n S_j) - f_n(\gamma^{-1}_n S_j + Z_j \varepsilon)) \\ \times \sum_{k=j+[j\delta]}^{n} \int f_{nkj} \, d(H^{\#}_{k-j} - H^{\#}_{k-j} * G_{\varepsilon_{nkj}}) \Bigg]$$

converges to 0, first by letting $n \to \infty$ and then $\varepsilon \to 0$. By taking $\upsilon = H^{\#}_{k-j}$, $\mu = H^{\#}_{k-j} * G_{\varepsilon_{nkj}}$ and $K_\eta = K_{\nu^{-1}_n \gamma_n \gamma^{-1}_{k-j}}$ in Lemma 14, one gets upper and lower bounds for

$$(28) \quad \int f_{nkj} \, d(H^{\#}_{k-j} * G_{\varepsilon_{nkj}} - H^{\#}_{k-j}).$$

By looking at these bounds, it is clear that we need to obtain bounds for the following:

$$(29) \int M_{f_{nkj}, \nu^{-1}_n \gamma_n \gamma^{-1}_{k-j}} \, d(H^{\#}_{k-j} * G_{\varepsilon_{nkj}} * K_{\nu^{-1}_n \gamma_n \gamma^{-1}_{k-j}} - H^{\#}_{k-j} * K_{\nu^{-1}_n \gamma_n \gamma^{-1}_{k-j}}),$$

$$(30) \int (M_{f_{nkj}, 2\nu^{-1}_n \gamma_n \gamma^{-1}_{k-j}} - f_{nkj}) \, dH^{\#}_{k-j},$$

$$(31) \int (f_{nkj} - m_{f_{nkj}, 2\nu^{-1}_n \gamma_n \gamma^{-1}_{k-j}}) \, dH^{\#}_{k-j}.$$

Regarding (30), we have, again using Lemma 14 with $\upsilon = 0$ and $K_\eta = K_{\nu^{-1}_n \gamma_n \gamma^{-1}_{k-j}}$,

$$\int (M_{f_{nkj}, 2\nu^{-1}_n \gamma_n \gamma^{-1}_{k-j}} - f_{nkj}) \, dH^{\#}_{k-j}$$

$$\leq \int (M_{f_{nkj},3\nu_n^{-1}\gamma_n\gamma_{k-j}^{-1}} - m_{f_{nkj},\nu_n^{-1}\gamma_n\gamma_{k-j}^{-1}})\, d(H_{k-j}^{\#} * K_{\nu_n^{-1}\gamma_n\gamma_{k-j}^{-1}})$$

$$(32) \qquad \leq \int (M_{f_{nkj},3\nu_n^{-1}\gamma_n\gamma_{k-j}^{-1}}(y) - m_{f_{nkj},3\nu_n^{-1}\gamma_n\gamma_{k-j}^{-1}}(y))\, dy$$

$$\times \int |\widehat{H}_{k-j}(u)||\widehat{K}_{\nu_n^{-1}\gamma_n\gamma_{k-j}^{-1}}(u)|\, du$$

$$\leq C\gamma_n\gamma_{k-j}^{-1}\left(\int (M_{f_n,\nu_n^{-1}}(y) - m_{f_n,\nu_n^{-1}}(y))\, dy\right)(1 + \rho^{k-j}\nu_n\gamma_n^{-1}\gamma_{k-j})$$

$$= C(\gamma_n\gamma_{k-j}^{-1} + \nu_n\rho^{k-j})\int (M_{f,\beta_n\nu_n^{-1}}(y) - m_{f,\beta_n\nu_n^{-1}}(y))\, dy,$$

where we have used (22), Lemma 15(i) $[S_{k-j}^{*}$ is in the form (15)$]$ and the identity (20). In the same way, the final bound (32) holds for (31) also. Now consider (29). Let $h_{k-j}(y)$ be the Lebesgue density of $H_{k-j}^{\#} * K_{\nu_n^{-1}\gamma_n\gamma_{k-j}^{-1}}$. Then

$$(33) \qquad \int M_{f_{nkj},\nu_n^{-1}\gamma_n\gamma_{k-j}^{-1}}\, d(H_{k-j}^{\#} * G_{\varepsilon_{nkj}} * K_{\nu_n^{-1}\gamma_n\gamma_{k-j}^{-1}} - H_{k-j}^{\#} * K_{\nu_n^{-1}\gamma_n\gamma_{k-j}^{-1}})$$

$$= \int M_{f_{nkj},\nu_n^{-1}\gamma_n\gamma_{k-j}^{-1}}(y)\left[\int h_{k-j}(y - z\varepsilon_{nkj})\phi(z)\, dz - h_{k-j}(y)\right]dy.$$

Now [recall that $H_{k-j}^{\#}(y) = H_{k-j}(y - \gamma_{k-j}^{-1}\gamma_n S_{nkj}^{**})$]

$$\left|\int h_{k-j}(y - z\varepsilon_{nkj})\phi(z)\, dz - h_{k-j}(y)\right|$$

$$= \left|\int (e^{-u^2\varepsilon_{nkj}^2/2} - 1)e^{iu(y - \gamma_{k-j}^{-1}\gamma_n S_{nkj}^{**})}(\widehat{H}_{k-j} * \widehat{K}_{\nu_n^{-1}\gamma_n\gamma_{k-j}^{-1}})(u)\, du\right|$$

$$\leq \int |e^{-u^2\varepsilon_{nkj}^2/2} - 1||(\widehat{H}_{k-j} * \widehat{K}_{\nu_n^{-1}\gamma_n\gamma_{k-j}^{-1}})(u)|\, du.$$

Using this and Lemma 15(ii), (33) is bounded in absolute value by

$$\left(\int M_{|f_{nkj}|,\nu_n^{-1}\gamma_n\gamma_{k-j}^{-1}}(y)\, dy\right)\left(\int |e^{-u^2\varepsilon_{nkj}^2/2} - 1||(\widehat{H}_{k-j} * \widehat{K}_{\nu_n^{-1}\gamma_n\gamma_{k-j}^{-1}})(u)|\, du\right)$$

$$(34) \leq C(\varepsilon_{nkj}^2 + \rho^{k-j}\nu_n\gamma_n^{-1}\gamma_{k-j})\gamma_n\gamma_{k-j}^{-1}\left(\int M_{|f|,\beta_n\nu_n^{-1}}(y)\, dy\right)$$

$$\leq C(\varepsilon_{nkj}^2\gamma_n\gamma_{k-j}^{-1} + \nu_n\rho^{k-j}).$$

Thus, summing (32) and (33), (28) is bounded in absolute value by a constant multiple of $(Q(n) + \varepsilon_{n,k-j}^2)\gamma_n\gamma_{k-j}^{-1} + (Q(n) + 1)\nu_n\rho^{k-j}$, where

$$Q(n) = \int (M_{f,\beta_n\nu_n^{-1}}(y) - m_{f,\beta_n\nu_n^{-1}}(y))\, dy.$$

In addition, $E[|f_n(\gamma_n^{-1} S_j) - f_n(\gamma_n^{-1} S_j + Z_j \varepsilon)|] \leq C(\gamma_n \gamma_j^{-1} + \nu_n \rho^j)$ by Proposition 19 and (23). Thus, (27) is bounded in absolute value by a constant multiple of

$$
(35) \quad \frac{1}{n^2} \sum_{j=[n\delta]}^{n} (\gamma_n \gamma_j^{-1} + \nu_n \rho^j)
$$
$$
\times \sum_{k=j+[j\delta]}^{n} ((Q(n) + \varepsilon_{nkj}^2)\gamma_n \gamma_{k-j}^{-1} + (Q(n)+1)\nu_n \rho^{k-j}).
$$

Now note that $\varepsilon_{nkj} = \varepsilon n^H (k-j)^{-H} u(n)(u(k-j))^{-1}$, where the slowly varying $u(n)$ is such that $\gamma_n = n^H u(n)$. When $k \geq j + [j\delta]$ and $j \geq [n\delta]$, there is a constant $C(\delta)$ depending on $\delta$ such that $n^H (k-j)^{-H} \leq C(\delta)$ for all sufficiently large $n$. Hence, one can check that $\frac{1}{n^2} \sum_{j=[n\delta]}^{n} \gamma_n \gamma_j^{-1} \sum_{k=j+[j\delta]}^{n} \varepsilon_{n,k-j}^2 \gamma_n \gamma_{k-j}^{-1} \leq C_*(\delta)\varepsilon^2$ for some constant $C_*(\delta)$. In addition, $\frac{1}{n^2} \sum_{j=1}^{n} \gamma_n \gamma_j^{-1} \sum_{k=j+1}^{n} \nu_n \rho^{k-j} \to 0$, $\frac{1}{n^2} \sum_{j=1}^{n} \nu_n \rho^j \sum_{k=j+1}^{n} \gamma_n \gamma_{k-j}^{-1} \to 0$ and $\frac{1}{n^2} \sum_{j=1}^{n} \nu_n \rho^j \sum_{k=j+1}^{n} \nu_n \rho^{k-j} \to 0$. (Recall that $n^{-1}\nu_n \to 0$.) Also, $Q(n) \to 0$ by (7). Thus, for each $\delta > 0$, (27) converges to 0 as $n \to \infty$ and $\varepsilon \to 0$. Hence, it remains to show that the difference between (26) and (27) converges to 0, first by letting $n \to \infty$ and then $\delta \to 0$. The same arguments used in Proposition 19 show that $|\int f_n(\gamma_n^{-1}\gamma_{k-j}y)\, d(H_{k-j}^\# * G_{\varepsilon_{nkj}} - H_{k-j}^\#)(y)|$ is bounded by $C(\gamma_n \gamma_{k-j}^{-1} + \nu_n \rho^{k-j})$. Hence, it can be seen that it is enough to show that $\frac{1}{n^2} \sum_{j=1}^{[n\delta]} \gamma_n \gamma_j^{-1} \sum_{k=j+1}^{n} \gamma_n \gamma_{k-j}^{-1}$ and $\frac{1}{n^2} \sum_{j=[n\delta]}^{n} \gamma_n \gamma_j^{-1} \times \sum_{k=j+1}^{j+[j\delta]} \gamma_n \gamma_{k-j}^{-1}$ converge to 0, first by letting $n \to \infty$ and then $\delta \to 0$. This is true for the first of these because it is bounded by $\frac{1}{n^2} \sum_{j=1}^{[n\delta]} \gamma_n \gamma_j^{-1} \sum_{l=1}^{n} \gamma_n \gamma_l^{-1}$. The same is true for the second one because it can be rewritten as

$$
\frac{1}{n^2} \sum_{j=[n\delta]}^{n} \gamma_n \gamma_j^{-1} \sum_{l=1}^{[j\delta]} \gamma_n \gamma_l^{-1} \sim C \frac{\gamma_n^2}{n^2} \sum_{j=[n\delta]}^{n} \gamma_j^{-1}[j\delta]\gamma_{[j\delta]}^{-1}
$$
$$
\sim C \frac{\delta^{1-H}\gamma_n^2}{n^2} \sum_{j=[n\delta]}^{n} j\gamma_j^{-2} \sim C\delta^{1-H} \qquad \text{as } n \to \infty,
$$

for each $\delta > 0$, where we have used the fact $j\gamma_j^{-2}$ is regularly varying with index $1 - 2H > -1$. This completes the proof of Proposition 6 under the assumptions of Theorem 2.

Now consider the proof under the assumptions of Theorem 3. First, consider Theorem 3(i). Because $\sup_y |f(y)| \leq C$, we have $\frac{1}{n} \sum_{j=1}^{2n_0-1} |f_n(\gamma_n^{-1} S_j + x)| \leq \frac{Cn_0\beta_n}{n}$, where $\frac{\beta_n}{n} \to 0$. Hence, one can restrict to the sum $\frac{1}{n} \sum_{j=2n_0}^{n} f_n(\gamma_n^{-1} S_j + x)$. Proceeding as in the proof of Proposition 19 and using Lemma 16(i) and the first part of Lemma 17, *for any* $r > 0$, there is an $L_r$ such that $E[|f_n(\gamma_n^{-1} S_k + x)|]$ is bounded by

$$\int |f_n(\gamma_n^{-1} \gamma_k y)| \, dy \left( \int |\hat{H}_k^*(u)| \, du \right) + \frac{C\beta_n L_r}{k^r}$$

$$\leq C\gamma_n \gamma_k^{-1} \int |f(y)| \, dy + \frac{C\beta_n L_r}{k^r}$$

for all $k \geq 2n_0$, where $H_k^*$ is as in Lemma 16 corresponding to the distribution $H_k$ of $S_k^*$. Hence, (24) (with the sums restricted to $2n_0 \leq k \leq n$) also holds, because $\int f^2(y) \, dy \leq (\sup_y |f(y)|) \int |f(y)| \, dy$. It remains to deal with (26), where as above, one can restrict the sum to $j$ and $k$ such that $k > j + 2n_0$ and $j \geq 2n_0$. In the same way as above but using the Lemma 17(ii), (28) is bounded by $C\gamma_n \gamma_{k-j}^{-1}(\varepsilon_{nkj}^2 + \rho^{k-j}) + \frac{C\beta_n L_r}{(k-j)^r}$ when $k > j + 2n_0$. Thus, the same arguments used earlier under the assumptions of Theorem 2 become applicable when $r > 1$. [Now $\gamma_n \gamma_{k-j}^{-1} \rho^{k-j} + \frac{\beta_n L_r}{(k-j)^r}$ plays the role of the earlier $\nu_n \rho^{k-j}$, where only the facts $n^{-1}\nu_n \to 0$ and $\sum_{k=j+1}^n \rho^{k-j} < C$ were used.] This completes the proof of Proposition 6 under Theorem 3(i).

Regarding the proof under Theorem 3(ii), note that we are considering the sum $\frac{1}{n} \sum_{j=n_0}^n f_n(\gamma_n^{-1} S_j + x)$ restricted to $n_0 \leq j \leq n$. Here we use the conclusions of Lemma 18, which hold for $H_j$ itself for $j \geq n_0$. The bound in Proposition 19 now becomes $C\gamma_n \gamma_k^{-1} \int |f(y)| \, dy$ and the conclusion (24) (with the sums restricted to $n_0 \leq k \leq n$) also holds. Also, in dealing with (26) one can restrict the sum to $j$ and $k$ such that $k > j + n_0$ (and $j \geq n_0$), because

$$|E[f_n(\gamma_n^{-1} S_j + x) f_n(\gamma_n^{-1} S_k + x)]|$$

$$\leq E^{1/2}[f_n^2(\gamma_n^{-1} S_j + x)] E^{1/2}[f_n^2(\gamma_n^{-1} S_k + x)]$$

so that, similar to (24), $\frac{1}{n^2} \sum_{j=n_0}^n |E[f_n(\gamma_n^{-1} S_j + x) f_n(\gamma_n^{-1} S_{j+q} + x)]| \to 0$ for each fixed $1 \leq q \leq n_0$. Further, as above, (28) is now bounded by $C\gamma_n \gamma_{k-j}^{-1}(\varepsilon_{nkj} + \rho^{k-j})$ when $k > j + n_0$. Hence, the proof of Proposition 6 is concluded. $\square$

REMARK 8. Under the conditions of Theorem 3(ii), it is implicit in the preceeding proof that, for every $0 < s_1 < s_2 \leq 1$, $\sup_x E[\{\frac{1}{n} \sum_{j=[ns_1]}^{[ns_2]} f_n(\gamma_n^{-1} \times S_j + x)\}^2]$ is bounded by $\frac{C}{n} \sum_{j=[ns_1]}^{[ns_2]} \gamma_n \gamma_j^{-1}(\frac{\beta_n}{n} \int f^2(y) \, dy + \frac{1}{n} \sum_{k=j+1}^{[ns_2]} \gamma_n \gamma_{k-j}^{-1})$. Similar bounds hold under Theorem 3(i) and under Theorem 2. One can establish analogous bounds for $\sup_x E[\{\frac{1}{n} \sum_{j=[ns_1]}^{[ns_2]} f_n(\gamma_n^{-1} S_j + x)\}^{2l}]$ for integers $l \geq 1$.

INDIAN STATISTICAL INSTITUTE
BANGALORE CENTER
8TH MILE, MYSORE ROAD
BANGALORE 560059
INDIA
E-MAIL: jegan@isibang.ac.in